\magnification=\magstep1

          %%%%%%%%%%%%%%%%%%%%%%%%%%%%%%%%%%%%%%%%%%%%%%%%%%%%%%%%%%%
          %       FOP.tex FORTEX-compatible version of papmacs      %
          %   VERSION OF May 6, 1989; does not use the AMS-fonts    %
          %%%%%%%%%%%%%%%%%%%%%%%%%%%%%%%%%%%%%%%%%%%%%%%%%%%%%%%%%%%

\def\item{\vskip1.3pt\hang\textindent}% THIS REPLACES KNUTH'S DEF'N

% THIS
                                              %REPLACES KNUTH'S DEF'N

\tolerance=300
\pretolerance=200
\hfuzz=1pt
\vfuzz=1pt

%\magnification=\magstep1
\hoffset 0cm            %all offsets went wrong
%\voffset=0.8true cm 
\hsize=5.8 true in
\vsize=9.5 true in

\def\rightheadline{\hfil\smc\lastname\hfil\tenbf\folio}
\def\leftheadline{\tenbf\folio\hfil\smc\lastname\hfil}
\headline={\ifodd\pageno\rightheadline\else\leftheadline\fi}
\newdimen\dimenone
\def\checkleftspace#1#2#3#4#5{%DIESER MACRO STAMMT VON APPELT
 \dimenone=\pagetotal
 \advance\dimenone by -\pageshrink   %testen ob Titel noch mit Gewalt auf Seite
                                                                          %geht
 \ifdim\dimenone>\pagegoal          %nacha tua nix-- gewoehnliche Outputroutine
   \else\dimenone=\pagetotal
        \advance\dimenone by \pagestretch
        \ifdim\dimenone<\pagegoal
          \dimenone=\pagetotal
          \advance\dimenone by#1         %addieren Skip vor Ueberschrift (=#1)
          \setbox0=\vbox{#2\parskip=0pt                %#2 ist gewaehlter Font
                       \hyphenpenalty=10000
                       \rightskip=0pt plus 5em
                       \noindent#3 \vskip#4}    %#3=Ueberschrift,#4=skip nachher
        \advance\dimenone by\ht0
        \advance\dimenone by 3\baselineskip
        \ifdim\dimenone>\pagegoal\vfill\eject\fi
          \else\eject\fi\fi}

\parindent=35pt
\mathsurround=1pt
\parskip=1pt plus .25pt minus .25pt
\normallineskiplimit=.99pt

\mathchardef\emptyset="001F % THIS REPLACES KNUTH'S DEFINITION

\def\Res{\mathop{\rm Res}\limits}
% USED FOR REAL PART OF COMPLEX NUMBERS
% USED FOR IMAGINARY PART OF COMPLEX NUMBERS
 % USED FOR IDENTITY FUNCTION
\def\Int{\mathop{\rm int}\nolimits}
% USED FOR TRACE OF MATRIX
%

%\def\u{\mathop{\rm u}\nolimits}

%\def\O{\mathop{\rm O}\nolimits}

\def\1{{\bf1}}\def\0{{\bf0}}

\def\({\bigl(}  \def\){\bigr)}
\def\<{\mathopen{\langle}}\def\>{\mathclose{\rangle}}

\def\Z{{\mathchoice{{\hbox{$\rm Z\hskip 0.26em\llap{\rm Z}$}}}%
{{\hbox{$\rm Z\hskip 0.26em\llap{\rm Z}$}}}%
{{\hbox{$\scriptstyle\rm Z\hskip 0.31em\llap{$\scriptstyle\rm Z$}$}}}{{%
\hbox{$\scriptscriptstyle\rm Z$\hskip0.18em\llap{$\scriptscriptstyle\rm Z$}}}}}}

\def\R{{\mathchoice{\hbox{$\rm I\hskip-0.14em R$}}%
{\hbox{$\rm I\hskip-0.14em R$}}%
{\hbox{$\scriptstyle\rm I\hskip-0.14em R$}}%
{\hbox{$\scriptscriptstyle\rm I\hskip-0.10em R$}}}}

\def\K{{\mathchoice{\hbox{$\rm I\hskip-0.15em K$}}%
{\hbox{$\rm I\hskip-0.15em K$}}%
{\hbox{$\scriptstyle\rm I\hskip-0.15em K$}}%
{\hbox{$\scriptscriptstyle\rm I\hskip-0.11em K$}}}}

\def\qed{\hfill {\hbox{[\hskip-0.05em ]}}}

\def\P{{\mathchoice{\hbox{$\rm I\hskip-0.14em P$}}%
{\hbox{$\rm I\hskip-0.14em P$}}%
{\hbox{$\scriptstyle\rm I\hskip-0.14em P$}}%
{\hbox{$\scriptscriptstyle\rm I\hskip-0.10em P$}}}}

\def\C{{\mathchoice%
{\hbox{$\rm C\hskip-0.47em\hbox{%
\vrule height 0.58em width 0.06em depth-0.035em}$}\;}%
{\hbox{$\rm C\hskip-0.47em\hbox{%
\vrule height 0.58em width 0.06em depth-0.035em}$}\;}%
{\hbox{$\scriptstyle\rm C\hskip-0.46em\hbox{%
$\scriptstyle\vrule height 0.365em width 0.05em depth-0.025em$}$}\;}%  \scripts
{\hbox{$\scriptscriptstyle\rm C\hskip-0.41em\hbox{%  is necessary for measures!
$\scriptscriptstyle\vrule height 0.285em width 0.04em depth-0.018em$}$}\;}}}

\def\Q{\QQ\,\,}

\def\.{{\cdot}}
\def\|{\Vert}
\def\ssk{\smallskip}
\def\msk{\medskip}
\def\bsk{\bigskip}
\def\giantskip{\vskip2\bigskipamount}

\def\giantbreak{\par \ifdim\lastskip<2\bigskipamount \removelastskip
         \penalty-400 \giantskip\fi}

\def\nin{\noindent}
\def\cen{\centerline}
\def\pagebreak{\vskip 0pt plus 0.0001fil\break}
\def\linebreak{\break}

\def\epsilon{\varepsilon}

\font\ninerm=cmr9
\font\eightrm=cmr8
\font\sixrm=cmr6

\font\eightbf=cmbx8
\font\sixbf=cmbx6

\font\eighti=cmmi8
\font\sixi=cmmi6
\font\ninesy=cmsy9
\font\eightsy=cmsy8
\font\sixsy=cmsy6

\font\eightit=cmti8

% SANS SERIF 10 POINT
 %SANS SERIF 10 POINT ITALIC

\font\eightsl=cmsl8

\font\eighttt=cmtt8
 %SLANTED TYPEWRITER 10 POINT
 %BOLD FACE MATH SYMBOLS 10 POINT
 %DUNHILL STYLE 10 POINT
 %SAN SERIF BOLD EXTENDED 10 POINT
 %USED FOR TITLES
 %USED FOR TITLES
\font\bfone=cmbx10 scaled\magstep1 %BOLDFACE AT MAGSTEP 1
 %BOLDFACE AT MAGSTEP 2
 %BOLDFACE AT MAGSTEP 3
\font\smc=cmcsc10
\font\itbig=cmmi10 scaled\magstep1

\font\small=cmcsc8

\def\no #1. {\bigbreak\vskip-\parskip\noindent\bf #1. \quad\rm}

\def\Proposition #1. {\checkleftspace{0pt}{\bf}{Theorem}{0pt}{}
\bigbreak\vskip-\parskip\noindent{\bf Proposition #1.}
\quad\it}

\def\Theorem #1. {\checkleftspace{0pt}{\bf}{Theorem}{0pt}{}
\bigbreak\vskip-\parskip\noindent{\bf  Theorem #1.}
\quad\it}
\def\Corollary #1. {\checkleftspace{0pt}{\bf}{Theorem}{0pt}{}
\bigbreak\vskip-\parskip\nin{\bf Corollary #1.}
\quad\it}
\def\Lemma #1. {\checkleftspace{0pt}{\bf}{Theorem}{0pt}{}
\bigbreak\vskip-\parskip\noindent{\bf  Lemma #1.}\quad\it}

\def\Definition #1. {\checkleftspace{0pt}{\bf}{Theorem}{0pt}{}
\rm\bigbreak\vskip-\parskip\noindent{\bf Definition #1.}
\quad}

\def\Remark #1. {\checkleftspace{0pt}{\bf}{Theorem}{0pt}{}
\rm\bigbreak\vskip-\parskip\noindent{\bf Remark #1.}\quad}

\def\Exercise #1. {\checkleftspace{0pt}{\bf}{Theorem}{0pt}{}
\rm\bigbreak\vskip-\parskip\noindent{\bf Exercise #1.}
\quad}

\def\Example #1. {\checkleftspace{0pt}{\bf}{Theorem}{0pt}{}
\rm\bigbreak\vskip-\parskip\noindent{\bf Example #1.}\quad}
\def\Examples #1. {\checkleftspace{0pt}{\bf}{Theorem}{0pt}
\rm\bigbreak\vskip-\parskip\noindent{\bf Examples #1.}\quad}

\newcount\problemnumb \problemnumb=0
\def\Problem{\global\advance\problemnumb by 1\bigbreak\vskip-\parskip\noindent
{\bf Problem \the\problemnumb.}\quad\rm }

\def\Proof#1.{\rm\par\ifdim\lastskip<\bigskipamount\removelastskip\fi\smallskip
            \noindent {\bf Proof.}\quad}

\nopagenumbers

\def\author{}
\def\lastname{}
\def\thanks#1{\footnote*{\eightrm#1}}
\def\title{}

\def\lastname{}
\def\h{{\textstyle{1\over2}}} 
 
\def\he{{1\over2}} 
\def\si{\sigma}
   
\def\n{{\cal N}}   
\def\ep{\epsilon}

\def\text{\textstyle} 
\def\disp{\displaystyle} 
\def\d{{\,\rm d}}

\def\and{{\rm and }}

\def\n{\cen{{\it W.G. Nowak}}}

\expandafter\edef\csname amssym.def\endcsname{%
       \catcode`\noexpand\@=\the\catcode`\@\space}
%  Set the catcode to 11 for use in private control sequence names.
\catcode`\@=11
\def\undefine#1{\let#1\undefined}
\def\newsymbol#1#2#3#4#5{\let\next@\relax
 \ifnum#2=\@ne\let\next@\msafam@\else
 \ifnum#2=\tw@\let\next@\msbfam@\fi\fi
 \mathchardef#1="#3\next@#4#5}
\def\mathhexbox@#1#2#3{\relax
 \ifmmode\mathpalette{}{\m@th\mathchar"#1#2#3}%
 \else\leavevmode\hbox{$\m@th\mathchar"#1#2#3$}\fi}
\def\hexnumber@#1{\ifcase#1 0\or 1\or 2\or 3\or 4\or 5\or 6\or 7\or 8\or
 9\or A\or B\or C\or D\or E\or F\fi}

%Loading fontfiles `eufm' and `msbm'
\font\tenmsb=msbm10
\font\sevenmsb=msbm7
\font\fivemsb=msbm5
\newfam\msbfam
\textfont\msbfam=\tenmsb
\scriptfont\msbfam=\sevenmsb
\scriptscriptfont\msbfam=\fivemsb
\edef\msbfam@{\hexnumber@\msbfam}
\def\Bbb#1{{\fam\msbfam\relax#1}}

\newsymbol\Bbbk 207C
\def\widehat#1{\setbox\z@\hbox{$\m@th#1$}%
 \ifdim\wd\z@>\tw@ em\mathaccent"0\msbfam@5B{#1}%
 \else\mathaccent"0362{#1}\fi}
\def\widetilde#1{\setbox\z@\hbox{$\m@th#1$}%
 \ifdim\wd\z@>\tw@ em\mathaccent"0\msbfam@5D{#1}%
 \else\mathaccent"0365{#1}\fi}
\font\teneufm=eufm10
\font\seveneufm=eufm7
\font\fiveeufm=eufm5
\newfam\eufmfam
\textfont\eufmfam=\teneufm
\scriptfont\eufmfam=\seveneufm
\scriptscriptfont\eufmfam=\fiveeufm

\catcode`@=11 %Set the catcode to 11 for use in private control sequence names.

\expandafter\edef\csname amssym.def\endcsname{%
       \catcode`\noexpand\@=\the\catcode`\@\space}
\font\eightmsb=msbm8
\font\sixmsb=msbm6
\font\fivemsb=msbm5
\font\eighteufm=eufm8
\font\sixeufm=eufm6
\font\fiveeufm=eufm5
\newskip\ttglue
\def\eightpoint{\def\rm{\fam0\eightrm}%
  \textfont0=\eightrm \scriptfont0=\sixrm \scriptscriptfont0=\fiverm
  \textfont1=\eighti \scriptfont1=\sixi \scriptscriptfont1=\fivei
  \textfont2=\eightsy \scriptfont2=\sixsy \scriptscriptfont2=\fivesy
  \textfont3=\tenex \scriptfont3=\tenex \scriptscriptfont3=\tenex
\textfont\eufmfam=\eighteufm
\scriptfont\eufmfam=\sixeufm
\scriptscriptfont\eufmfam=\fiveeufm
\textfont\msbfam=\eightmsb
\scriptfont\msbfam=\sixmsb
\scriptscriptfont\msbfam=\fivemsb
  \def\it{\fam\itfam\eightit}%
  \textfont\itfam=\eightit
  \def\sl{\fam\slfam\eightsl}%
  \textfont\slfam=\eightsl
  \def\bf{\fam\bffam\eightbf}%
  \textfont\bffam=\eightbf \scriptfont\bffam=\sixbf
   \scriptscriptfont\bffam=\fivebf
  \def\tt{\fam\ttfam\eighttt}%
  \textfont\ttfam=\eighttt
  \tt \ttglue=.5em plus.25em minus.15em
  \normalbaselineskip=9pt
  \def\MF{{\manual opqr}\-{\manual stuq}}%
  \let\big=\eightbig
  \setbox\strutbox=\hbox{\vrule height7pt depth2pt width\z@}%
  \normalbaselines\rm}
\def\eightbig#1{{\hbox{$\textfont0=\ninerm\textfont2=\ninesy
  \left#1\vbox to6.5pt{}\right.\n@space$}}}

%  Restore the catcode value for @ that was previously saved.

\csname amssym.def\endcsname

% end of \symb1.tex

\def\la{\lambda} 
\def\al{\alpha}

\def\om{\omega} 
\def\({\left(} 
\def\){\right)} 
\def\for{\qquad \hbox{for}\ } 
\def\eq{\eqalign} 
\def\f{{1\over 2\pi i}} 
 
\def\O#1{O\(#1\)} 
\def\abs#1{\left| #1 \right|}

\def\klein{\eightpoint \def\smc{\small} \baselineskip=9pt}   

\def\fn#1#2{{\parindent=0.7true cm
\footnote{$^{(#1)}$}{{\klein  #2}}}}

\font\boldmas=msbm10                  %%
\def\Bbb#1{\hbox{\boldmas #1}}        %%
\def\Z{{\Bbb Z}}                        %%
\def\Q{{\Bbb Q}}
\def\R{{\Bbb R}}

\def\C{{\Bbb C}}

                  %%
        %%

%%%%%%%%%%%%%%%%%%%%%%%%%%%%%%%%%%%%%%%%%%%%%%%%%%%%%%%%%%%%%%%%%%%%%%
\font\eightrm=cmr8                                                    
\long\def\fussnote#1#2{{\baselineskip=9pt                            
\setbox\strutbox=\hbox{\vrule height 7pt depth 2pt width 0pt}% 
\eightrm                                                         
\footnote{#1}{#2}}}                                              
%%%%%%%%%%%%%%%%%%%%%%%%%%%%%%%%%%%%%%%%
%% zum Verkleinern (Summationsgrenzen, Folgenindex, e.a.)%%
%%%%%%%%%%%%%%%%%%%%%%%%%%%%%%%%%%%%%%%%
\font\boldmasi=msbm10 scaled 700      %%
\def\Bbbi#1{\hbox{\boldmasi #1}}      %%
\font\boldmas=msbm10                  %%
\def\Bbb#1{\hbox{\boldmas #1}}        %%
\def\Pi{{\Bbbi P}}                      %%
\def\Qi{{\Bbbi Q}}                      %%

\def\Ki{{\Bbbi K}}

                      %%
%%%%%%%%%%%%%%%%%%%%%%%%%%%%%%%%%%%%%%%%

                        %%

\def\dint #1 {% Doppelintegral; #1 Integrationsbereich 
\quad  \setbox0=\hbox{$\disp\int\!\!\!\int$}
  \setbox1=\hbox{$\!\!\!_{#1}$}
  \vtop{\hsize=\wd1\centerline{\copy0}\copy1} \quad}

\def\drint #1 {% Dreifachintegral; #1 Integrationsbereich 
\qquad  \setbox0=\hbox{$\disp\int\!\!\!\int\!\!\!\int$}
  \setbox1=\hbox{$\!\!\!_{#1}$}
  \vtop{\hsize=\wd1\centerline{\copy0}\copy1}\qquad}

\def\frac#1#2{{#1\over #2}}

\def\date{\the\day.~\the\month.~\the\year}

\def\mod{\,{\rm mod}\,}  
\def\klein{\eightpoint \def\smc{\small} }

\def\at#1#2#3{{\left. #1 \right|}_{#2}^{#3}}

\def\frac#1#2{{#1\over#2}} 
\def\Int{\int\limits}

%\nonumbers
\hsize=16true cm     \vsize=23.7true cm

\parindent=0cm

\def\E{{\cal A}} 
\def\A{{\cal A}}
\def\r{{\cal R}}
\def\p{{\cal P}}
\def\I{{\cal I}}
\def\n{{\cal N}}
\def\c{{\cal C}}
\def\M{{\cal M}}
\def\f{{1\over2\pi i}} 
\def\B{{\cal B}}
\def\ok{{\cal O}_\Ki} 
\def\Cat{{\bf C}} 
\def\cat{L(2,\chi_D)}

\vbox{\vskip 1.5true cm}   

\font\bigma=cmr10 scaled \magstep1

\cen{{\bfone The average number of solutions}}\ssk 
\cen{{\bfone \textfont0=\bigma 
of the Diophantine equation $\hbox{\itbig U}^{\hbox{\itbig 2}} + 
\hbox{\itbig V}^{\,\hbox{\itbig 2}} 
= \hbox{\itbig W}^{\,\hbox{\rm 3}}$}}\msk  
\cen{{\bfone and related arithmetic functions}}  \bsk \msk 

\cen{\bf M.~K\"uhleitner and W.G.~Nowak} 

\vbox{\vskip 1.5true cm}   
{\klein {\it Abstract. } 
For the number of integer solutions of the title equation, with 
$W\le x$ ($x$ a large parameter), an asymptotics of the form $A x\log x + B x 
+\O{x^{1/2}(\log x)^3(\log\log x)^2}$ is established. This is achieved in a 
general setting which furnishes applications to some other natural arithmetic 
functions.\msk 
{\it AMS-Classification. } 11N37, 11M06, 11R42, 11D25}

\vbox{\vskip 1.5true cm}   

{\bf 1.~Introduction. }  During the problem session of the 1991 Czechoslovak 
Number Theory Conference, A.~Schinzel proposed the following question: Let 
$r(n)$ denote the number of ways to write the positive integer $n$ as a sum of 
two squares, and consider the asymptotic formula 
$$ \sum_{n\le x} (r(n))^2 = 4 x\log x + C x + \O{x^{1/2}\om(x)}\,,\eqno(1.1) $$ 
as $x\to\infty$. How small can the factor $\om(x)$ be made, if one uses the 
sharpest tools of contemporary analytic number theory? \ssk 
In fact, asymptotics of the shape\fn{1}{To complete the history, one should 
mention the older, somewhat coarser results of W.~Sierpi\'nski [16] and 
S.~Ramanujan [12], with error terms $\O{x^{3/4}\log x}$, resp., 
$\O{x^{3/5+\ep}}$. Furthermore, Schinzel [15] himself had bounded the 
remainder from below, showing that it is an $\Omega(x^{3/8})$.} 
(1.1) had been established by B.M.~Wilson 
[19], with $\om(x)=x^\ep$, and W.~Recknagel [13], with $\om(x)=(\log x)^6$. 
\par 
Subsequently, as a reaction to Schinzel's proposal, the first named author [6] 
sharpened the estimate to $\om(x)=(\log x)^{11/3}(\log\log x)^{1/3}$. \ssk 
The problem addressed by Schinzel is closely related to another arithmetic 
question to which our present title refers: Given a large parameter $x$, how 
many integer triples $(u,v,w)$ exist with $u^2+v^2=w^3$ , $w\le x$ ? This 
matter was dealt with by K.H.~Fischer [2] and also by Recknagel [13]. \ssk 
Since ${1\over4}r(\cdot)$ is multiplicative, it is easy to write up the 
corresponding generating Dirichlet series for $(r(n))^2$ and $r(n^3)$. (See 
section 4.2 for details.) Thus, 
both problems are subsumed in a natural way by the following more general 
result which is to be the objective of the present note. \bsk

{\bf Theorem.}\quad{\it Let $a(n)$ be an arithmetic function satisfying 
$a(n)\ll n^{\ep} $ for every $\ep>0$, with a Dirichlet series 
$$  F(s) = \sum_{n=1}^\infty {a(n)\over n^s} = {(\zeta_\Ki(s))^2 \over 
(\zeta(2s))^{m_1} (\zeta_\Ki(2s))^{m_2}}\,G(s) \qquad (\Re(s)>1)\,,  $$ 
where $\zeta_\Ki$ is the Dedekind zeta-function of some quadratic number 
field $\K$, $G(s)$ is holomorphic and bounded in some half-plane 
$\Re(s)\ge\theta$, $\theta<\h$, and $m_1, m_2$ are nonnegative integers. 
Then, for $x$ large, 
$$ \eq{\sum_{n\le x} a(n) &= \Res_{s=1}\(F(s)\, {x^s \over s}\) 
+ \O{x^{1/2}(\log x)^3 (\log\log x)^{m_1+m_2}} = \cr 
&= A x \log x + B x + \O{x^{1/2}(\log x)^3 (\log\log x)^{m_1+m_2}}\,. \cr} $$}

\msk  

{\bf Remark. } It is natural to ask why we have chosen just this very degree 
of generality, resp., specialization, in our suppositions. Note that, for 
$\K$a quadratic field, $\zeta_\Ki(s)=\zeta(s)L(s,\chi)$, where $\chi$ is a 
certain real Dirichlet character. Thus the above function $F$ is a special 
case of 
$$ {L(s,\chi^{(1)}) \dots L(s,\chi^{(J)}) \over \zeta_{\Ki_1}(2s) \dots 
\zeta_{\Ki_M}(2s)}\,G(s)\,, $$ 
where the $L(s,\chi^{(j)})$ are ordinary Dirichlet $L$-functions, 
$\K_1,\dots,\K_M$ are arbitrary algebraic number fields, and $G(s)$ is as 
before. As long as $J=4$, our 
argument would go through without major alterations and lead to the same 
result, apart from the exponents of the log- and loglog-factors in the error 
term. \ssk 
However, for $J\ne4$, the situation changes drastically: If $J = 2$ or $3$, 
an important r\^ole is played by the zero-free region of the 
denominator.\fn{2}{Therefore, in these cases, 
the assumption of the truth of the Riemann Hypothesis (RH) leads to a better 
error term $O(x^\theta)$, with some $\theta<\h$. In contrast, our Theorem is 
as sharp as it would be if RH could be proven.} 
Using the sharpest information available of this kind, one obtains a 
bound $O(x^{1/2}\exp(-c(\log x)^{3/5}(\log\log x)^{-1/5}))$ for the remainder, 
where the exp-factor is familiar from the prime number theorem. (See, for 
instance, formulas (14.28) -- (14.30) and (1.105) -- (1.107) in Ivi\'c [5].) 
\par Opposed to 
this, for $J>4$ the size of the nominator in vertical strips becomes 
important. This ultimately leads to an estimate $O(x^{\al_J+\ep})$, where 
$\al_J$ is the best known error-exponent in the $J$-dimensional Piltz divisor 
problem which, for $J\ge5$, is still a matter of small improvements from time 
to time (cf.~Ivi\'c [5], p.~355, or Titchmarsh [18], ch.~12). 
Thus the case $J=4$ remains as the most delicate one, at least as far as 
smaller factors, apart from $x^{1/2}$, in the error term are concerned. \ssk 
In addition, 
our analysis has been inspired by recent work\fn{3}{The authors are indebted 
to Professor Sankaranarayanan for sending them a preprint of this paper.} 
of K.~Ramachandra and A.~Sankaranarayanan [11] which deals with 

\vbox{$$ \sum_{n=1}^\infty {(d(n))^2\over n^s} = {(\zeta(s))^4\over\zeta(2s)} 
\qquad (\Re(s)>1)\,.  $$ 
See also Sankaranarayanan's Oberwolfach lecture [14]. However, our argument is 
technically a bit simpler. }

\bsk\msk

{\bf 2.~Some auxiliary results.} \bsk 

{\bf Lemma.}\quad {\it Let $\K$ be an arbitrary number field of degree 
$[\K:\Q]\ge1$, $\zeta_\Ki$ its Dedekind zeta-function, and $\ep>0$ fixed. 
Then, for each sufficiently large $T$ there exists a measurable set 
$\E_T=\E_{T,\Ki}\subset[T,2T]$ of Lebesgue measure 
$$  \la(\E_T) \le T^\ep\,,  \eqno(2.1)  $$ 
with the property that, for some $C_{\ep,\Ki}>0$ depending only on $\ep$ and 
$\K$, 
$$  \sup_{t\in[T,2T]\setminus\E_T}\abs{\zeta_\Ki(1+it)}^{-1} \le C_{\ep,\Ki}\,
\log\log T\,. \eqno(2.2)  $$  } \bsk 
{\bf Proof. } For the case of the Riemann zeta-function, this is contained in 
Ramachandra [10], Theorem 1. Our argument follows the lines of Lemma 3.2 in 
Ramachandra and Sankaranarayanan [11]. \ssk According to Heath-Brown [3], there 
exists some $\delta=\delta(\ep)>0$ such that, for $T$ large, the cardinality 
of the set 
$$  \r_{\ep,T} = \{\rho\in\C :\ \zeta_\Ki(\rho)=0,\ \Re(\rho)\ge1-3\delta,\  
T\le\Im(\rho)\le2T\ \}  \eqno(2.3)  $$ 
satisfies $$\#\r_{\ep,T} \ll T^{\ep/2}\,. \eqno(2.4)  $$
For any fixed $c>0$, we define 
$$  \A(c,T) := [T, T+cT^{\ep/4}] \cup [2T-cT^{\ep/4}, 2T] \cup 
\bigcup_{\rho\in\r_{\ep,T}} [\Im(\rho)- cT^{\ep/4}, \Im(\rho)+cT^{\ep/4}]\,,$$ 
and 
$$ \M(c,T) := \{s\in\C : \Re(s)\ge1-3\delta,\,
\Im(s)\in[T,2T]\setminus\A(c,T)\, \} \cup \{s\in\C : \Re(s)\ge1\, \}\,. 
\eqno(2.5) $$ 
On every $\M(c,T)$, by construction $\zeta_\Ki(s)\ne0$, thus 
$\log\zeta_\Ki(s)$ can be defined properly by analytic continuation\fn{4}{To 
be explicite, on $\M(c,T)$, 
we can define $\log\zeta_\Ki(s) := \log\zeta_\Ki(2)  
+ \int_\c (\zeta_\Ki'(z)/\zeta_\Ki(z)) \d z$ where $\c$ consists of the two 
straight line segments from 2 to $2+\Im(s)i$ and from $2+\Im(s)i$ to $s$.}. 
For any $\eta=2+it$ with $t\in[T,2T]\setminus\A(1,T)$, 
we consider the circular discs 
$$ \c_1 := \{s\in\C:\ \abs{s-\eta}\le1+\delta\ \}\,,\quad 
\c_2 := \{s\in\C:\ \abs{s-\eta}\le1+2\delta\ \}\,.  $$ 
Evidently, $\c_1\subset\c_2\subset\M(\h,T)$. Thus 
the Borel-Carath\'eodory inequality (cf.~Titchmarsh [17], ch.~5.5) yields 
$$  \max_{s\in\c_1}\abs{\log\zeta_\Ki(s)} \le 
{2(1+\delta)\over\delta}\max_{s\in\c_2} \log\abs{\zeta_\Ki(s)} 
+ {2+3\delta\over\delta}\abs{\log\zeta_\Ki(2+it)}\,.  $$ 
Since, for a certain $C>0$, \ $\abs{\zeta_\Ki(\si+it)}\le\abs{t}^C$, 
provided that $\si\ge\h$, $\abs{t}$ sufficiently large, it follows that 
$\log\abs{\zeta_\Ki(s)} \le C \log(2T)$ on $\c_2$. Hence, in particular, 
$$  \abs{\log\zeta_\Ki(\si+it)} \ll \log T \for \si\ge1-\delta\,, $$ 
the constant involved not depending\fn{5}{Note that, by (2.7), 
$\abs{\log\zeta_\Ki(2+it)}\le\log\zeta_\Ki(2)$.} on $t$. 
Therefore,  
$$  \log\zeta_\Ki(s) \ll \log T \eqno(2.6)  $$ 
uniformly for all $s$ with $\Re(s)\ge1-\delta$, 
$\Im(s)\in[T,2T]\setminus\A(1,T)$. \ssk 
We now use the series representation (valid for $\Re(s)>1$) 
$$  \log\zeta_\Ki(s) = \sum_{\p}\sum_{m=1}^\infty {1\over m} \n(\p)^{-ms} = 
\sum_{\p} \n(\p)^{-s} + H(s)\,, \eqno(2.7)  $$ 
where $H(s)$ is regular and bounded in any half-plane $\Re(s)\ge\theta>\h$. 
($\p$ are the prime ideals in the ring $\ok$ of algebraic integers in the 
field $\K$ and $\n(\cdot)$ denotes the norm of ideals in $\ok$.) \ssk 
We put $X=(\log T)^{1/\delta}$. Then, for any $t\in\R$, 
$$ \f\Int_{\delta-i\infty}^{\delta+i\infty} \log\zeta_\Ki(1+it+w) 
\Gamma(w) X^w \d w = $$ 
$$ = \sum_{\p}\sum_{m=1}^\infty {1\over m} \n(\p)^{-m(1+it)} 
\(\f\Int_{\delta-i\infty}^{\delta+i\infty} \(X^{-1}\n(\p)^m\)^{-w} 
\Gamma(w)\d w\) = $$ 
$$ = \sum_{\p}\sum_{m=1}^\infty {1\over m} \n(\p)^{-m(1+it)}\exp(-\n(\p)^m/X) = 
\sum_{\p} \n(\p)^{-1-it}\exp(-\n(\p)/X) + O(1)\,. \eqno(2.8) $$ 
To estimate the remaining sum, we observe that 
$$ \sum_{\n(\p)> X} \n(\p)^{-1-it}\exp(-\n(\p)/X) \ll 
\sum_{\n(\I)> X} \n(\I)^{-1} \exp(-\n(\I)/X) = $$ 
$$  = \Int_X^\infty {1\over u} e^{-u/X}\,\d\(\sum_{\n(\I)\le u} 1\)\ 
\ll 1\,, $$ using only that $\disp\sum_{\n(\I)\le u} 1 \ll u$. ($\I$ denotes 
arbitrary integral ideals in $\ok$.) Further, 
$$ \abs{\sum_{\n(\p)\le X} \n(\p)^{-1-it}\exp(-\n(\p)/X)} \le 
\sum_{\n(\p)\le X} \n(\p)^{-1} = \log\log X + O(1)\,,  $$ 
where the last conclusion is immediate from the prime ideal theorem in the 
form (cf.~Narkiewicz [8], pp.~369-372) 
$$  \sum_{\n(\p)\le u} 1 = {u\over\log u} + \O{{u\over(\log u)^2}}\,. $$ 
In view of (2.8) and the choice $X=(\log T)^{1/\delta}$, this implies that 
$$ \abs{\f\Int_{\delta-i\infty}^{\delta+i\infty} \log\zeta_\Ki(1+it+w) 
\Gamma(w) X^w \d w} \le \log\log\log T + O(1)\,.  \eqno(2.9)  $$ 
Our next step is to put $W = (\log\log T)^2$ and to recall Stirling's formula 
in the weak form (valid uniformly in any strip $\si_1\le\si\le\si_2$, 
$\abs{t}\ge1$)
$$\Gamma(\si+it) \ll \abs{t}^{\si-1/2} \exp(-{\text{\pi\over2}}\abs{t})\,. 
\eqno(2.10)  $$ 
From this it readily follows that  
$$ \Int_{\delta\pm iW}^{\delta\pm i\infty} \log\zeta_\Ki(1+it+w) 
\Gamma(w) X^w \d w \ll X^\delta \Int_W^\infty 
\exp(-{\text{\pi\over2}} u) \d u \ll X^\delta 
\exp(-{\text{\pi\over2}} W) \ll 1\,, $$   
hence 
$$ \abs{\f\Int_{\delta-iW}^{\delta+iW} \log\zeta_\Ki(1+it+w) 
\Gamma(w) X^w \d w} \le \log\log\log T + O(1)\,.  \eqno(2.11)  $$  
From now on we impose the condition that $t\in[T,2T]\setminus\A(2,T)$. 
We evaluate the integral in (2.11) by the residue theorem, 
applied to the rectangle $\Re$ with vertices $\pm\delta\pm iW$. (Note that 
for $t\in[T,2T]\setminus\A(2,T)$ and $w\in\Re$, 
necessarily $1+it+w\in\M(1,T)$, and for $\log(\zeta_\Ki(1+it+w))$ 
the bound (2.6) applies.) We obtain 
$$  \Int_{-\delta-iW}^{-\delta+iW} \log\zeta_\Ki(1+it+w) 
\Gamma(w) X^w \d w \ll X^{-\delta}\,\log T \Int_{-W}^W 
\abs{\Gamma(-\delta+iu)}\d u \ll 1\,,  $$ 
and similarly, for the horizontal segments, 
$$ \Int_{-\delta\pm iW}^{\delta\pm iW} \log\zeta_\Ki(1+it+w) 
\Gamma(w) X^w \d w \ll X^{\delta}\,\log T 
\exp(-{\text{\pi\over2}} W) \ll 1\,. $$   
Since the only pole inside the rectangle, at $w=0$, gives a residue of 
$\log\zeta_\Ki(1+it)$, we altogether derive from (2.11) that \par 
\vbox{$$  \abs{\log\zeta_\Ki(1+it)} \le \log\log\log T + C_1\,,  $$ 
for all $t\in[T,2T]\setminus\A(2,T)$ ($C_1>0$ an appropriate constant 
depending on $\K$ and $\ep$).} \par Hence, for the same $t$, 
$$ \abs{\zeta_\Ki(1+it)}^{-1} = \exp(-\log\abs{\zeta_\Ki(1+it)}) \le 
\exp(\abs{\log\zeta_\Ki(1+it)}) \le e^{C_1}\,\log\log T\,.  $$ 
Taking $\A_{T,\Ki} = \A(2,T)$, this just proves clause (2.2) of the Lemma, 
while (2.1) is immediate by (2.4) and (2.5). \qed \msk  
We conclude this section by stating some more bounds\fn{6}{In fact, 
the estimates (2.12) and (2.13) are far away from being the sharpest ones of 
their kind, but they suffice for our purpose and are conveniently available in 
textbooks.} for the zeta-functions 
involved, which will be needed in the proof of the Theorem. First of all, 
again for every number field $\K$ with $[\K:\Q]\ge1$, 
$$  \abs{\zeta_\Ki(\si+it)}^{-1} \ll (\log(2+\abs{t}))^C \eqno(2.12)  $$ 
uniformly (at least) in $\si\ge1$, the constant $C>0$ possibly depending on 
$\K$. (This is most conveniently deduced after the classic example of Apostol 
[1], Th.~13.7, using the necessary facts about the Dedekind zeta-function 
from Narkiewicz [8], ch.~7.2.) \ssk 
Further, we recall that 
$$ \zeta(\h+it) \ll (1+\abs{t})^{1/6+\ep}\,,\qquad 
L(-\ep'+it,\chi) \ll_\chi\ (1+\abs{t})^{1/2+\ep} $$ 
for any Dirichlet $L$-series. (For the first estimate, see Titchmarsh 
[18], Th.~5.12. The second one follows readily from the 
functional equation in Apostol [1], Th.~12.11, along with Stirling's formula.) 
Applying the Phragm\'en-Lindel\"of theorem in the form given by 
Titchmarsh [17], ch.~5.65, we infer that 
$$ \zeta(\si+it) \ll 1+\abs{t}^{(1-\si)/3+\ep}\,, \qquad 
L(\si+it,\chi)\ll_\chi\ 1+\abs{t}^{(1-\si)/2+\ep}\,,  $$  
uniformly (at least) in $\si\ge\h$, $\abs{t}\ge1$. 
Now, if $\K$ is a {\it quadratic } field of discriminant $D$, 
its zeta-function can be factorized as $\zeta_\Ki(s) = \zeta(s) L(s,\chi_D)$. 
(Cf.~Zagier [21], p.~100, in particular for the definition of the real 
character $\chi_D$ involved.) Therefore, in this case, 
$$ \zeta_\Ki(\si+it) \ll 1+\abs{t}^{{5\over6}(1-\si)+\ep}\,, \eqno(2.13) $$ 
uniformly in $\si\ge\h$, $\abs{t}\ge1$. 
Last but not least, it has been proved by W.~M\"uller 
[7] (even in the stronger form of an asymptotics) that 
$$  \Int_0^T \abs{\zeta_\Ki(\h+it)}^2 \d t \ll T(\log T)^2\,, \eqno(2.14) $$ 
again for every quadratic field $\K$. 

\bsk\msk

\vsize=23.2true cm

{\bf 3.~Proof of the Theorem.}\quad W.l.o.g., let $x$ be half an odd integer. 
Then, by a simple version of Perron's truncated formula (see, e.g., Prachar 
[9], p.~376), 
$$ \sum_{n\le x} a(n) = \f \Int_{1+\ep-ix^{3/5}}^{1+\ep+ix^{3/5}} 
F(s)\, x^s {\d s\over s} + \O{x^{1/2}}\,,  \eqno(3.1)  $$ 
for any fixed $\ep>0$ sufficiently small. We shift the path of integration to 
the vertical line $\Re(s)=\h$. The horizontal segments contribute 
$$ \f\Int_{\he\pm i x^{3/5}}^{1+\ep\pm i x^{3/5}} F(s)\, x^s {\d s\over s} 
\ll \Int_\he^{1+\ep} 
{\abs{\zeta_\Ki(\si\pm i x^{3/5})}^2\over\abs{\zeta(2\si\pm 2i x^{3/5})}^{m_1} 
\abs{\zeta_\Ki(2\si\pm 2i x^{3/5})}^{m_2}}\,x^{\si-3/5}\,\d\si $$ 
By (2.12) and (2.13), this is 
$$ \ll x^{-3/5+\ep} \max_{\he\le\si\le1+\ep} \(x^\si (x^{3/5})^{5(1-\si)/3} \) 
\ll x^{2/5+\ep}\,.  $$ 
Hence, by the residue theorem, 
$$ \sum_{n\le x} a(n) = \Res_{s=1}\(F(s)\, {x^s \over s}\) + 
\f \Int_{\he-ix^{3/5}}^{\he+ix^{3/5}} 
F(s)\, x^s {\d s\over s} + \O{x^{1/2}}\,.  \eqno(3.2)  $$ 
To estimate the remaining integral, we use the parametrization $s=\h(1+it)$, 
$\abs{t}\le 2x^{3/5}$, along with a dyadic decomposition, to obtain 
$$  \Int_{\he-ix^{3/5}}^{\he+ix^{3/5}} F(s)\, x^s {\d s\over s} 
\ll x^{1/2}\(1 + \sum_{T=2^{-j}x^{3/5},\atop j=0,1,2,\dots} {1\over T} 
\Int_T^{2T} {\abs{\zeta_\Ki(\h(1+ it))}^2\over\abs{\zeta(1+it)}^{m_1} 
\abs{\zeta_\Ki(1+it)}^{m_2}}\,\d t \)\,,  \eqno(3.3)  $$ \ssk 
where the sum in fact contains only $O(\log x)$ terms. With $\A_{\Ki,T}$ as in 
the Lemma, 
we put $\B_{\Ki,T}=[T,2T]\setminus\A_{\Ki,T}$. By (2.1), (2.12), and (2.13),  
$$  \Int_{\disp\A_{\Ki,T}} 
{\abs{\zeta_\Ki(\h(1+ it))}^2\over\abs{\zeta(1+it)}^{m_1} 
\abs{\zeta_\Ki(1+it)}^{m_2}}\,\d t \ll T^{5/6+4\ep}\,.  $$ 
Moreover, by (2.2) and (2.14), 
$$  \Int_{\disp\B_{\Ki,T}} 
{\abs{\zeta_\Ki(\h(1+ it))}^2\over\abs{\zeta(1+it)}^{m_1} 
\abs{\zeta_\Ki(1+it)}^{m_2}}\,\d t \ll T (\log T)^2 (\log\log T)^{m_1+m_2}\,. 
$$ Inserting the last two bounds into (3.3) and summing over $T$, we obtain 
$$  \Int_{\he-ix^{3/5}}^{\he+ix^{3/5}} F(s)\, x^s {\d s\over s} \ll x^{1/2} 
(\log x)^3 (\log\log x)^{m_1+m_2}\,.  $$ 
Together with (3.2), this completes the proof of our Theorem. \bsk\msk 

{\bf 4.~Applications.}\quad It remains to verify that the result established 
covers the two special problems addressed in the title and in the 
introduction, even in a more general context. \ssk 
{\bf 4.1.~The second moment of quadratic Dedekind-zeta coefficients. } 
For an arbitrary quadratic number field $\K$ with discriminant $D$, 
let $\ok$ denote the ring of algebraic integers in $\K$, 
and $r_\Ki(n)$ the number of integral ideals $\I$ in $\ok$ of Norm 
$\n(\I)=n$. Then, as we shall show below, for $\Re(s)>1$, 
$$  \sum_{n=1}^\infty {(r_\Ki(n))^2\over n^s} = 
{(\zeta_\Ki(s))^2\over\zeta(2s)}\, \prod_{p|D}(1+p^{-s})^{-1}\,.  
\eqno(4.1)  $$ 
Therefore, applying the Theorem with $(m_1,m_2)=(1,0)$, 
we obtain what follows\fn{7}{For the evaluation of the residue at $s=1$ and 
the subsequent numerical computations, we have 
employed {\it Mathematica } [20].}. \bsk 
{\bf Corollary 1.}\quad{\it For any quadratic field $\K$ of discriminant $D$, 
and $x$ large, 
$$\sum_{n\le x} (r_\Ki(n))^2 
= A_1 x \log x + B_1 x + \O{x^{1/2}(\log x)^3 \log\log x}\,, $$
with 
$$ A_1 =  {6\over\pi^2} L(1,\chi_D)^2 \prod_{p|D} {p\over p+1}\,, $$ 
$$ B_1 =  {6\over\pi^2} L(1,\chi_D)^2 \prod_{p|D} {p\over p+1}\,
\(-1 +2\gamma + \sum_{p|D} {\log p\over p+1} \ + {2 L'(1,\chi_D)\over 
L(1,\chi_D)} - {12\over\pi^2} \zeta'(2)\)\,,  $$ 
where $\gamma$ is the Euler-Mascheroni constant.} \msk 
Since $r(n) = 4 r_{\Qi(i)}(n)$, this implies in particular that 
$$ \sum_{n\le x} (r(n))^2 
= 4 x \log x + 16 B_1 x + \O{x^{1/2}(\log x)^3 \log\log x} \eqno(4.2) $$ 

\vbox{with $$ 16 B_1 = -4+8\gamma +{4\over3}\log2 +{32\over\pi} L'(1) - 
{48\over\pi^2}\zeta'(2) \approx 8.0665\,,$$ 
where $L(\cdot)$ is the 
$L$-function corresponding to the non-principal Dirichlet character mod 4.} 
\ssk 
It remains to verify (4.1). To this end, we recall the decomposition laws in 
$\ok$ (cf.~Zagier [21]). On the set $\P$ of all rational primes $p$, there 
exists a partition 
$\P=\P_0\cup\P_1\cup\P_2$ such that\fn{8}{Note that, on the basis of 
this partition, the Dirichlet character $\chi_D$ can be defined as follows: 
$\chi(p)=0$ if $p\in\P_0$, $\chi(p)=1$ if $p\in\P_1$, and 
$\chi(p)=-1$ if $p\in\P_2$.} \ssk 
\vbox{$$ p\in\P_0 \iff p|D \iff (p)=\p^2\,,\ \n(\p)=p\,, \leqno(i)$$ 
$$ p\in\P_1 \iff (p) = \p_1 \p_2\,,\ \p_1\ne\p_2\,,\ \n(\p_1)=\n(\p_2)=p\,, 
\leqno(ii) $$ 
$$ p\in\P_2 \iff (p) \hbox{ prime in }\ok\,,\ \n(p)=p^2\,.\leqno(iii)  $$ }

For $\Re(s)>1$, this implies the Euler product representation 
$$  \eq{\zeta_\Ki(s) &= \prod_{p\in\Pi} 
\(1+\sum_{j=1}^\infty r_\Ki(p^j)p^{-js}\) 
= \prod_{\p} (1-\n(\p)^{-s})^{-1} = \cr 
&= \prod_{p\in\P_0}(1-p^{-s})^{-1} \prod_{p\in\P_1}(1-p^{-s})^{-2} 
\prod_{p\in\P_2}(1-p^{-2s})^{-1}\,. \cr  }  \eqno(4.3)  $$ 
By comparison, we immediately see that $r_\Ki(p^j)=1$ throughout for 
$p\in\P_0$, \ $r_\Ki(p^j)=j+1$ for $p\in\P_1$, and $r_\Ki(p^j)=\cases{1 & if 
$j$ is even,\cr 0 & if $j$ is odd,\cr}$\quad for $p\in\P_2$. \ssk 
Hence, 
$$  \sum_{n=1}^\infty {(r_\Ki(n))^2\over n^s} = 
\prod_{p\in\P_0}(1-p^{-s})^{-1} \prod_{p\in\P_1}\(1+\sum_{j=1}^\infty (j+1)^2 
p^{-js}\) \prod_{p\in\P_2}(1-p^{-2s})^{-1}\,.  $$ 
Using that \ 
$\disp\sum_{j\ge0}(j+1)^2 z^j = (1+z)(1-z)^{-3}$ (for $\abs{z}<1$), 
we see that this equals 

\vbox{$$ \prod_{p\in\P_0}(1-p^{-s})^{-1} 
\prod_{p\in\P_1} (1-p^{-s})^{-3} (1+p^{-s})  
\prod_{p\in\P_2}(1-p^{-2s})^{-1} = $$ 
\cen{$ = {\disp\prod_{p\in\Pi}}(1-p^{-2s}) 
\({\disp\prod_{p\in\P_0}}(1-p^{-s})^{-1} {\disp\prod_{p\in\P_1}} 
(1-p^{-s})^{-2} 
{\disp\prod_{p\in\P_2}}(1-p^{-2s})^{-1}\)^2 
{\disp\prod_{p\in\Pi_0}}(1+p^{-s})^{-1} = $}   
$$  =  {(\zeta_\Ki(s))^2\over\zeta(2s)}\, \prod_{p|D}(1+p^{-s})^{-1}\,, $$ 
which proves (4.1).} \msk 
{\bf 4.2.~Diophantine equations as addressed in the title. } Again for $\K$ a 
quadratic field, we shall show that 
$$  F(s) := \sum_{n=1}^\infty {r_\Ki(n^3)\over n^s} = 
{(\zeta_\Ki(s))^2\over\zeta(2s)\zeta_\Ki(2s)}\,G(s) \quad 
(\Re(s)>1)\,,  \eqno(4.4) $$ 
where $G(s)$ is holomorphic and bounded in every half-plane 
$\Re(s)\ge\theta>{1\over3}$. Taking this for granted, we apply our Theorem 
with $(m_1,m_2)=(1,1)$ and derive the following consequence. \bsk 
{\bf Corollary 2.}\quad {\it For any quadratic field $\K$ of discriminant $D$, 
and $x$ large, 
$$\sum_{n\le x} r_\Ki(n^3)
= A_2 x \log x + B_2 x + \O{x^{1/2}(\log x)^3 (\log\log x)^2}\,, $$
with 
$$ A_2 = {36\,L(1,\chi_D)^2\,G(1) \over\pi^4\,\cat}\,, $$ 
$$ B_2 = {36\,L(1,\chi_D)^2\,G(1)\over\pi^4\, \cat}\,
\({2 L'(1,\chi_D)\over L(1,\chi_D)}- 1 +
2\gamma +{G'(1)\over G(1)}-{2 L'(2,\chi_D)\over\cat}-
{24\zeta'(2)\over\pi^2}\)\,, $$                 
where $G(s)$ is given in form of an Euler product in $(4.6)$ 
below}\fn{9}{Of course, $G'(1)/G(1)$ is most easily evaluated as 
$\at{{\d\over\d s}\log G(s)}{s=1}{}$, which transforms the products into sums 
of derivatives of logarithms.}. \bsk 

To apply this result to Diophantine equations, let $Q=Q(u,v)=au^2+buv+cv^2$ be 
an integral, primitive, positive definite binary quadratic form of class 
number 1.\fn{10}{Equivalently, the discriminant $D=b^2-4ac$ is one of the 
"{\it numeri idonei}" 
$\{-3,-4,-7,-8,-11,$ $ -12, -16, -19, -27,-28,-43,-67,-163\}$. 
See, e.g., Hlawka/Schoiáengeier [4], p.~92.} Then it is well-known 
(cf.~Zagier [21], \S\S\ 8 and 11) that 
$$ r_Q(n) := \#\{(u,v)\in\Z^2:\ Q(u,v)=n\ \} 
= \om_D\,r_\Ki(n)\,,  \eqno(4.5)  $$ 
where $\K=\Q(\sqrt{D})$, and\fn{11}{$\om_D$ is the number of automorphisms of 
the form $Q$ or, in other terms, the number of units in 
${\cal O}_{\Qi(\sqrt{D})}$.} $\om_D=\cases{6 & for $D=-3$, \cr 4 & for 
$D=-4$, \cr 2 & else. \cr }$ \msk 
We thus can infer the following conclusion. \bsk 
\vbox{{\bf Corollary 3.}\quad{\it For every 
integral, primitive, positive definite binary quadratic form $Q$ of class
number 1, with discriminant $D$, and for large real $x$, 
$$ \eq{&\#\{(u,v,w)\in\Z^3:\ Q(u,v) = w^3,\ w\le x\ \} = \cr 
=&\ 1+\sum_{1\le n\le x} r_Q(n^3) = 1+ \om_D\,\sum_{1\le n\le x} 
r_{\Qi(\sqrt{D})}(n^3) =\cr 
=&\ \om_D(A_2 x \log x + B_2 x) + \O{x^{1/2}(\log x)^3 (\log\log x)^2}\,, \cr} 
$$ where $A_2$, $B_2$ are as in Corollary $2$.}} \msk 
In particular, for the Diophantine equation of the title, we obtain 
$$ \eq{&\#\{(u,v,w)\in\Z^3:\ u^2+v^2 = w^3,\ w\le x\ \} =  
 1+\sum_{1\le n\le x} r(n^3) =\cr 
&\ = 4A_2 x \log x + 4B_2 x + \O{x^{1/2}(\log x)^3 (\log\log x)^2}\,, \cr} $$ 
with 
$$ 4A_2 = {9 G(1)\over\pi^2\, \Cat}\approx 0.9091 \,, $$ 
$$ 4B_2 = {9 G(1)\over\pi^2\, \Cat}\,\({8L'(1,\chi_{-4})\over\pi}- 1 +
2\gamma + {G'(1)\over G(1)}-{2 L'(2,\chi_{-4})\over\Cat}-
{24\zeta'(2)\over\pi^2}\) \approx 2.1715  \,, $$ 
where $\disp\Cat=L(2,\chi_{-4})=\sum_{j\ge0}(-1)^j (2j+1)^{-2}$ is Catalan's 
constant, and, by (4.6) below, 
$$ G(1) = {8\over9}\, \prod_{p\equiv 1 \mod 4} {p^3(p+2)\over(p-1)(p+1)^3}\,
\prod_{p\equiv 3 \mod 4} {p^4\over p^4-1} \approx 0.91317\,, $$  
$$ {G'(1)\over G(1)} = -{1\over3} \log2 
- 6\,\sum_{p\equiv 1 \mod 4} {\log p\over(p+1)(p-1)(p+2)} 
- 4 \,\sum_{p\equiv 3 \mod 4} {\log p\over p^4-1} \approx -0.35876 \,.  $$ \bsk

It remains to verify (4.4). By what we noted earlier (right after (4.3)), for 
$\Re(s)>1$, 
$$  F(s) = \sum_{n=1}^\infty {r_\Ki(n^3)\over n^s} = 
\prod_{p\in\Pi} \(1 + \sum_{j=1}^\infty r_\Ki(p^{3j}) p^{-js}\) = $$ 
$$ = \prod_{p\in\P_0}(1-p^{-s})^{-1} \prod_{p\in\P_1} 
\(1+\sum_{j=1}^\infty (3j+1) p^{-js}\) \prod_{p\in\P_2}(1-p^{-2s})^{-1}\,. $$ 
Using that \ $\disp\sum_{j\ge0} (3j+1) z^j = (1+2z)(1-z)^{-2}$ 
(for $\abs{z}<1$), we may write this as 
$$ F(s) = \prod_{p\in\P_0}(1-p^{-s})^{-1} \prod_{p\in\P_1} 
{1+2p^{-s}\over(1-p^{-s})^2} \prod_{p\in\P_2}(1-p^{-2s})^{-1}\,. $$ 

\vbox{Recalling the Euler product representation of $\zeta_\Ki(s)$ and 
$\zeta_\Ki(2s)$, as given in (4.3), along with that of $\zeta(2s)$, we obtain 
after a brief calculation 
$$ G(s) = F(s)\,{\zeta(2s)\zeta_\Ki(2s)\over(\zeta_\Ki(s))^2} = 
\prod_{p\in\P_0}{1-p^{-s}\over(1-p^{-2s})^2}  
\prod_{p\in\P_1} h(p^{-s})\, 
\prod_{p\in\P_2} (1-p^{-4s})^{-1}\,,  \eqno(4.6) $$ 
where $h(z):=(1+2z)(1-z)^2(1-z^2)^{-3} = 1 + O(\abs{z}^3)$ around $z=0$. Thus 
$G(s)$ possesses a Dirichlet series absolutely convergent for 
$\Re(s)>{1\over3}$. This proves (4.4).} \bsk 

%\vbox{\vskip 1.5true cm}   

\klein \parindent=0pt 

\cen{\bf References}  \bsk

[1] {\smc T.M.~Apostol,} Introduction to analytic number theory. 
New York, Springer, 1976. \ssk 

[2] {\smc K.H.~Fischer,} \"Uber die Anzahl der Gitterpunkte auf Kreisen mit 
quadratfreien Radienquadraten. Arch.~Math. {\bf 33} (1979), 150--154. \ssk 

[3] {\smc D.R.~Heath-Brown,} On the density of the zeros of the Dedekind 
zeta-function. Acta Arithm. {\bf 33} (1977), 169--181.  \ssk 

[4] {\smc E.~Hlawka \& J.~Schoi\ss engeier,} Zahlentheorie. 2.~Aufl., 
Wien, Manz, 1990.  \ssk 

[5] {\smc A.~Ivi\'c,} The Riemann zeta-function. New York, Wiley, 1985. \ssk 

[6] {\smc M.~K\"uhleitner,} On a question of A.~Schinzel concerning the sum 
$\disp\sum_{n\le x}(r(n))^2$. Proc.~\"Osterr.-Ung.-Slowak. Koll. 
Zahlentheorie, Graz-Mariatrost 1992. Grazer Math.~Ber. {\bf 318} (1993), 
63--67. \ssk 

[7] {\smc W.~M\"uller,} The mean square of the Dedekind zeta-function in 
quadratic number fields. Math.~Proc. Cambridge Phil.~Soc. {\bf 106} (1989), 
403--417.\ssk 

[8] {\smc W.~Narkiewicz,} Elementary and analytic theory 
of algebraic numbers. 2nd ed., Warszawa, Springer \& PWN (Polish Scientific 
Publishers), 1990.  \ssk 

[9] {\smc K.~Prachar,} Primzahlverteilung. Berlin, Springer, 1957. \ssk 

[10] {\smc K.~Ramachandra,} A large value theorem for $\zeta(s)$. 
Hardy-Ramanujan J. {\bf 18} (1995), 1--9.  \ssk 

[11] {\smc K.~Ramachandra \& A.~Sankaranarayanan,} On an asymptotic formula 
of Srinivasa Ramanujan. Acta Arithm., in press.  \ssk 

[12] {\smc S.~Ramanujan,} Some formulae in the analytic theory of numbers. 
Mess.~Math. {\bf 45} (1916), 81--84.  \ssk 

[13] {\smc W.~Recknagel,} Varianten des Gau\ss schen Kreisproblems. 
Abh.~Math.~Sem.~Univ.~Hamburg {\bf 59} (1989), 183--189.  \ssk 

[14] {\smc A.~Sankaranarayanan,} On an asymptotic formula of Srinivas 
Ramanujan. Conf.~"Theory of the Riemann zeta and allied functions", held at 
Oberwolfach, Sept.~2001, Math.~FI Oberwolfach, Report No.~43/2001, p.~15. 
\ssk 

[15] {\smc A.~Schinzel,} On an analytic problem considered by Sierpi\'nski and 
Ramanujan. Proc.~Int.~Conf.~in honour of V.~Kubilius, Palanga (Lithuania) 
1991, VSP BV Int.~Corp., 1992, pp.~165--171.  \ssk 

[16] {\smc W.~Sierpi\'nski,} Sur la sommation de la s\'erie 
$\disp\sum_{a<n\le b} \tau(n) f(n)$, o\'u $\tau(n)$ signifie le nombre de 
d\'ecompositions du nombre 
$n$ en une somme de deux carr\'es de nombres entiers. (In Polish, French 
summary.) Prace Mat.~Fiz. {\bf 18} (1908), 1--59.  \ssk 

[17] {\smc E.C.~Titchmarsh,} The theory of functions. 2nd ed., 
Oxford, Univ.~Press, 1939.  \ssk 

[18] {\smc E.C.~Titchmarsh,} The theory of the Riemann zeta-function. 
2nd ed., revised by D.R.~Heath-Brown. Oxford, Univ.~Press, 1986.  \ssk 

[19] {\smc B.M.~Wilson,} Proofs of some formulae enunciated by Ramanujan. 
Proc.~London Math.~Soc. {\bf 21} (1922), 235--255. \ssk 

[20] {\smc Wolfram} Research, Inc., Mathematica 4.1. Champaign 2001. \ssk 

[21] {\smc D.B.~Zagier,} Zetafunktionen und quadratische K\"orper. Berlin, 
Springer, 1981.  

\vbox{\vskip 1.5true cm}  

\parindent=1true cm

\vbox{Manfred K\"uhleitner \& Werner Georg Nowak

Institut f\"ur Mathematik u.~Ang.Stat. 

Universit\"at f\"ur Bodenkultur 

Peter Jordan-Stra\ss e 82 

A-1190 Wien, Austria \ssk 

E-mail: {\tt kleitner@edv1.boku.ac.at,  nowak@mail.boku.ac.at} \ssk 

Web: http://www.boku.ac.at/math/nth.html}

\bye